\documentclass{article}
\usepackage{graphicx} 
\usepackage{indentfirst}
\usepackage{amsmath}
\usepackage{geometry}
\usepackage{chemmacros}

\title{Process Flowsheet Optimization with Surrogate and Implicit Formulations of a Gibbs Reactor}
\author{
$^{1,2}$Sergio I. Bugosen,
$^2$Carl D. Laird,
$^{1,3,*}$Robert B. Parker
\\
\\
{\small$^1$\it Information Systems and Modeling, Los Alamos National Laboratory, Los Alamos, NM}\\
{\small$^2$\it Department of Chemical Engineering, Carnegie Mellon University, Pittsburgh, PA}\\
{\small$^3$\it Center for Nonlinear Studies, Los Alamos National Laboratory, Los Alamos, NM}\\
{\small$^*$\it Corresponding author}
}
\date{April 2024}

\begin{document}
\maketitle

\section{Abstract}

Alternative formulations for the optimization of chemical process flowsheets are presented that leverage surrogate models and implicit functions to replace and remove, respectively, the algebraic equations that describe a difficult-to-converge Gibbs reactor unit operation. Convergence reliability, solve time, and solution quality of an optimization problem are compared among full-space, ALAMO surrogate, neural network surrogate, and implicit function formulations. Both surrogate and implicit formulations lead to better convergence reliability, with low sensitivity to process parameters. The surrogate formulations are faster at the cost of minor solution error, while the implicit formulation provides exact solutions with similar solve time. In a parameter sweep on the autothermal reformer flowsheet optimization problem, the full-space formulation solves 33 out of 64 instances, while the implicit function formulation solves 52 out of 64 instances, the ALAMO polynomial formulation solves 64 out of 64 instances, and the neural network formulation solves 48 out of 64 instances. This work demonstrates the trade-off between accuracy and solve time that exists in current methods for improving convergence reliability of chemical process flowsheet optimization problems.


\section{Introduction}

The selection of optimal operating conditions is a fundamental task in chemical process design. This requires the minimization or maximization of an objective function subject to nonlinear and nonconvex constraints. While local solvers such as CONOPT \cite{Drud1985} and IPOPT \cite{Wchter2005} can handle these problems, their convergence can be very sensitive to model formulation, initial guess, and scaling factors. In addition, flowsheet design equations usually contain complex nonlinear expressions, including logarithms, high-degree polynomials and bilinear terms. At certain variable values, these equations can be undefined, and their Jacobian can become singular, hindering convergence. 

Some methods to improve convergence of a full space flowsheet optimization problem and address sensitivity to initial parameters include using (1) sophisticated initialization routines and (2) model reformulation strategies, such as surrogate models and implicit functions.

In the case of (1), it is well known that an effective initialization strategy determines how easily the optimization algorithm converges to a solution \cite{Mazzei2014}. However, finding a good initialization is a cumbersome task because some unit operations do not easily converge for a set of specifications, and there is no systematic way to determine good initial values \cite{Caballero2008}. In addition, there is a significant computational cost associated with trying multiple initialization methods for each problem instance one attempts to solve.

Regarding (2), surrogate models aim to be simple models that approximate the input and output behaviour of complex systems (in our case study, a Gibbs reactor) over a specific input domain. A surrogate model typically relaxes accuracy in exchange for lower dimensionality and more reliable convergence \cite{Thombre2015}. Even if the development of these models requires computationally expensive simulations over a range of input values, they are computationally cheaper to evaluate once they are embedded into a larger interconnected system of equations, such as a chemical process flowsheet. Two surrogates that have been recently used in chemical engineering are neural networks and ALAMO polynomials \cite{Bhosekar2018}. Neural networks have been widely used in process control, modeling, and optimization. For instance, Henao and Maravelias use neural networks to model the production of maleic anhydride in a superstructure optimization problem \cite{Henao2010}. The ALAMO framework \cite{Cozad2014} is a recently developed tool used to build simple and accurate polynomial surrogates from a minimal set of training data. It makes use of an integer programming technique to choose basis functions of the input variables and compute the output variables as a linear combination thereof \cite{Bhosekar2018, Ma2022}. Surrogate optimization using the ALAMO framework has been used in superstructure optimizations for carbon capture systems \cite{Miller2014}, in global optimization of polygeneration systems \cite{Subramanian2022}, and distillation sequences \cite{Ma2022}.

Reformulating a model with the implicit function theorem, as proposed by \cite{Parker2022}, aims to exploit the fact that the difficulty of converging a large-scale nonlinear program (NLP) may be due to the effort required to converge the system of nonlinear equality constraints corresponding to specific unit operations in the flowsheet. Solving these units separately from the original formulation can significantly improve convergence reliability. In this work, we take advantage of non-singularity of the Gibbs reactor equations. Given that the Jacobian of these equations with respect to the Gibbs reactor's output variables is non-singular, solving the Gibbs reactor equations yields a unique solution for these output variables. This implies that there exists a function mapping state and input variables to the outputs of the Gibbs reactor. The utility of this property is that complicated algebraic equations from a unit operation can be removed from the NLP using implicit functions. This leaves a smaller set of equations to be seen by the solver and fewer variables to initialize and scale \cite{Parker2022}.



The surrogate and implicit reformulation strategies are similar in that they both remove complicating equations that describe a difficult-to-converge unit model. In this work, we design an autothermal reformer flowsheet (ATR) in the IDAES process modeling framework \cite{Lee2021} and we formulate a full space optimization problem, as shown in Section 3. We then compare convergence reliability, solution quality, and solve time among full-space, implicit function, and surrogate formulations of this problem, and study the trade-offs among these formulations.

\section{Background}

The general structure of the three different formulations is displayed in this section.


\subsection{Full Space Formulation}

Vector $x$ is a vector of state variables corresponding to each unit model in the flowsheet, while $u$ is a vector of manipulated inputs. Function $G$ describes operational constraints, function $H$ describes the equality constraints corresponding to each unit model in the flowsheet not including the unit model we intend to replace, and function $R$ describes the equality constraints only for this unit model. Vector $y$ is a vector containing the outlet variables corresponding to this unit and any additional intermediate variables used only by this unit. Eq. (\ref{eq1}d) will be replaced by either a surrogate model or an implicit function.

\begin{subequations}\label{eq1}
  \begin{align}
    \max &~~ f(x, u) \\
    \text{s.t.} &~~ G(x, u) \leq 0 \\
    &~~ H(x, y, u) = 0 \\
    &~~ R(x, y, u) = 0
  \end{align}
\end{subequations}

\subsection{Surrogate Formulation}

This formulation is identical to the one shown in Eq. (\ref{eq1}) except for Eq. (\ref{eq2}d), which instead of representing first-principles design equations for the unit model of interest, contains a surrogate model $\overline{R}(x,y,u)$ that approximates this unit’s behavior. While this surrogate can be obtained using a variety of methods, as described in \cite{Bhosekar2018}, in this paper we use neural network and polynomial surrogates.

\begin{subequations}\label{eq2}
  \begin{align}
    \max &~~ f(x, u) \\
    \text{s.t.} &~~ G(x, u) \leq 0 \\
    &~~ H(x, y, u) = 0 \\
    &~~ \overline{R}(x,y,u) = 0
  \end{align}
\end{subequations}

\subsection{Implicit Formulation} 

The implicit function theorem states that Eq. (\ref{eq1}d) can be reformulated as $y=R_y(x,u)$ if $\nabla_y R$ is non-singular for all values of state and input variables. Under this condition, Eq. (\ref{eq3}) is an exact reformulation of Eq. (\ref{eq1}). We solve for $y$ externally as a square system of equations in a separate interface and the resulting values and derivatives are communicated back to the optimization solver, specifically to Eq. (\ref{eq3-equality}), which links the outlet variable $y$ calculated by the implicit function to the inner optimization problem. Since this formulation keeps Eq. (\ref{eq1}d) feasible and the unit model equations are not seen by the NLP solver, we expect it to lead to better convergence reliability.

\begin{subequations}\label{eq3}
  \begin{align}
    \max &~~ f(x, u) \\
    \text{s.t.} &~~ G(x, u) \leq 0 \\
    &~~ H(x, R_y(x, u), u) = 0\label{eq3-equality}
  \end{align}
\end{subequations}

To solve this formulation with a second-order optimization method (e.g., IPOPT), the routines implemented in \cite{Parker2022} are used to calculate the constraint Jacobian, objective gradient, and Hessian of the Lagrangian. As an example, Eq. (\ref{eq4}) shows the Jacobian of $y$ as a function of $x$ and $u$, which is used to calculate the former derivative matrices. 

\begin{equation}\label{eq4}
\nabla_{x,u}y = -\nabla_y R^{-1} \nabla_{x,u}R
\end{equation}

\section{Problem Statement}

\begin{figure}[h]
\centering
\resizebox{0.65\columnwidth}{!}{\includegraphics{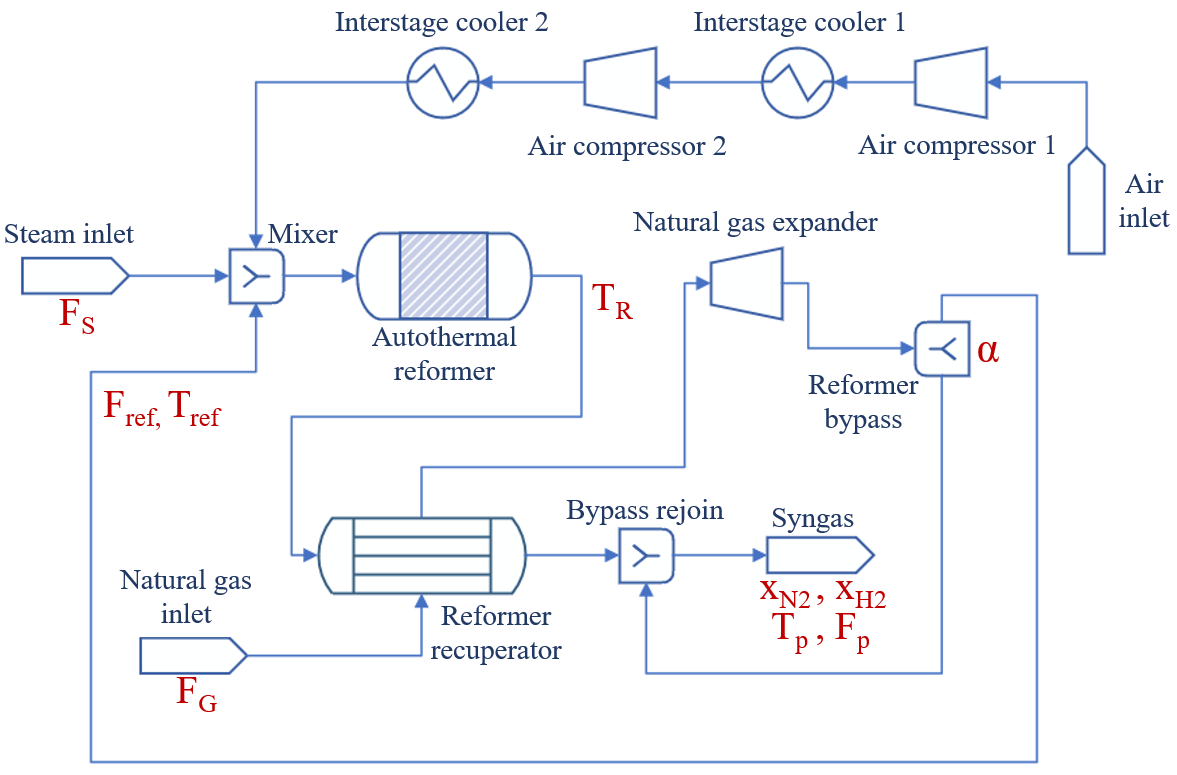}}
\caption{\label{fig:figure 1} Process Flow Diagram of the ATR process.}
\end{figure}

An autothermal reforming flowsheet is used as an example in the Optimization \& Machine Learning Toolkit (OMLT) \cite{OMLT}. The main objective of the autothermal reformer (ATR) is to produce syngas, mainly composed of \ch{H2}, \ch{CO}, \ch{CH4} and \ch{CO2}. The process is shown in Figure \ref{fig:figure 1}. First, a mixture of natural gas, steam and air is fed into the ATR (modeled as a Gibbs reactor). The hot syngas is then circulated through a shell and tube heat exchanger, also called the reformer recuperator, to heat the natural gas feed. This natural gas feed is then expanded to generate electrical power and is finally fed into the reactor, closing the loop. 

In this study, the objective of all optimization problems is to maximize the hydrogen concentration in the syngas stream, as shown in Eq. (\ref{eq6}),

\begin{equation}
  \max ~~ x_{\ch{H2}}\label{eq6}
\end{equation}

The operational constraints are specified in Eq. (\ref{eq7}). 

\begin{subequations}\label{eq7}
\begin{equation}
  T_{R} \leq 1200~[{\rm K}]
\end{equation}

\begin{equation}
  T_{p} \leq 650~[{\rm K}]
\end{equation}

\begin{equation}
F_{p} \geq 3500~\left[\frac{\rm mol}{\rm s}\right]
\end{equation}

\begin{equation}
  x_{\ch{N2}} \leq 0.3
\end{equation}

\begin{equation}
1120 \leq F_{G} \leq 1250~\left[\frac{\rm mol}{\rm s}\right]
\end{equation}

\begin{equation}
0 \leq \alpha \leq 1
\end{equation}

\begin{equation}
200 \leq F_{s} \leq 400~\left[\frac{\rm mol}{\rm s}\right]
\end{equation}

\end{subequations}

Where $T_{R}$ is the autothermal reformer outlet temperature and $T_{p}$, $F_{p}$ and $x_{\ch{N2}}$ correspond to product's temperature, molar flow rate and nitrogen concentration respectively. The manipulated variables are $F_{G}$, $F_{s}$ and $\alpha$, corresponding to inlet molar flow rate of natural gas, inlet molar flow rate of steam, and bypass fraction respectively. The full space optimization problem is composed of 898 variables and 895 equality constraints describing material, energy, and momentum balances for each unit model, where 55 of these correspond to algebraic equations describing the Gibbs reactor, such as Gibbs energy minimization, and enthalpy, element, pressure, total flow, and component flow balances. 

As an example, the Gibbs minimization equation is shown in Eq. (\ref{eq8}a). 

\begin{subequations}\label{eq8}
\begin{equation}
g_{partial,j} + \sum_e (L_{e}  \times \beta_{j,e}) = 0, ~\forall~j\in{J}
\end{equation}

\begin{equation}
J = \left\{\ch{H2}, \ch{CO}, \ch{H2O}, \ch{CO2}, \ch{CH4}, \ch{C2H6}, \ch{C3H8}, \ch{C4H10}, \ch{O2}\right\}
\end{equation}

\begin{equation}
J_{\rm inert} = \left\{\ch{N2}, \ch{Ar}\right\}
\end{equation}
\end{subequations}

Here, $g_{partial,j}$ is the partial molar Gibbs energy of component $j$, $L_{e}$ is the Lagrange multiplier of element $e$ and $\beta_{j,e}$ is the number of moles of element $e$ in one mole of component $j$ \cite{Lee2021}.

In this work, we implement full-space, surrogate, and implicit function formulations for solving this optimization problem. The Gibbs reactor, Eq. (\ref{eq1}d), will be replaced by a surrogate and an implicit function.


\section{Implementation} 

We used Pyomo 6.7.1 \cite{Bynum2021-xv}, an open-source optimization modeling language, and IDAES 2.4.0 \cite{Lee2021}, a process modeling framework, to design the ATR flowsheets in Python 3.9.6. The optimization problems were solved with IPOPT 3.14.11 \cite{Wchter2005} using linear solver MA27, called via the CyIpopt interface. The surrogates were obtained using the ALAMO machine learning framework \cite{Cozad2014} and TensorFlow 2.15. The neural network was embedded into the optimization problem using OMLT 1.1 \cite{OMLT}. Derivative computations are performed via the AMPL solver library (ASL) \cite{Gay1997} via the PyNumero interface \cite{Rodriguez2023}. Results were produced on a machine with an Apple M1 Max processor and 32 GB of RAM running macOS Ventura 13.6.6. 

\subsection{Dataset Generation}

A dataset describing the Gibbs reactor was generated to train the ALAMO polynomials and the neural network. This data was generated from 625 samples of the four-dimensional input space: five samples in each input arranged in a regular grid. The ranges of sampled inputs are shown in Table \ref{table:1}. The outputs (13 in total) correspond to reactor’s outlet temperature (Tout), outlet molar flow rate (Fout), and outlet compositions for the eleven components. Time to acquire this data was 400 seconds

\begin{table}[ht]
\centering
\begin{tabular}{ |c|c|c| } 
 \hline
 \textbf{Input variable} & \textbf{Unit} & \textbf{Range} \\
 \hline
Inlet steam molar flow rate ($F_S$) & [$\frac{mol}{s}$] & 200 - 350 \\ 
 Inlet natural gas temperature from reformer bypass ($T_{ref}$) & [$K$] & 600 - 900 \\
 Inlet natural gas flow rate from reformer bypass ($F_{ref}$) & [$\frac{mol}{s}$] & 600 - 900 \\
 \ch{CH4} Conversion in autothermal reformer ($X$)& [\%] & 80 - 95 \\
 \hline
\end{tabular}
\caption{Inputs and ranges of the ALAMO and neural network surrogates.}
\label{table:1}
\end{table}

This data set was partitioned into 80\% training and 20\% validation data to train the ALAMO and neural network surrogates and then gauge their accuracy. We consider that a surrogate is accurate if the coefficient of determination ($R^2$) is greater than 0.8 for each of the 13 parity plots obtained from the validation dataset. 

\subsection{ALAMO surrogate}

The Gibbs reactor was replaced by a surrogate block containing simple algebraic equations determined by the ALAMO machine learning framework. To balance the bias-variance trade-off and calculate a model with low nonlinear complexity, the four basis models we consider are a quadratic, a cubic, a linear variable and a constant. We have the option to include bilinear terms and higher degree polynomials to achieve a higher surrogate accuracy. However, for the purpose of this research, convergence could be hindered with the inclusion of those terms, particularly if we observe that a linear combination of simple basis functions can effectively approximate the reformer’s behavior. The surrogate model is com-posed of 13 equations and 3 variables, in contrast to the 55 equations and variables that model this first-principles Gibbs reactor. A subset of these 13 equations is displayed in Eq. (\ref{eq9}). Training time to acquire this surrogate model was 1.6 s.

 
\begin{subequations}\label{eq9}
\begin{equation}
T_{out} = 8.2 \times 10^{-4}F_S + 0.41 F_{ref}^3 + 897.4 X
\end{equation}

\begin{equation}
F_{out} = 3.9F_{ref} + 1.1F_S - 7.9 \times 10^{-7}F_{ref}^2 + 685X^2
\end{equation}

\begin{equation}
x_{\ch{H2}} = 5.3 \times 10^{-4} F_{ref} - 1.5\times 10^{-10} F_{ref}^3 + 0.14X^3
\end{equation}

\begin{equation}
x_{\ch{CO}} = -6.2 \times 10^{-10}F_S^3 + 8.1 \times 10^{-4} F_{ref} - 4.1 \times 10^{-7} F_{ref}^2 + 0.2X - 0.35
\end{equation}

\begin{equation}
x_{\ch{CH4}} = 1.5 \times 10^{-5}F_{ref} - 6 \times 10^{-6}F_S - 0.33X^2 + 0.16X^3 + 0.16
\end{equation}

\end{subequations}

\subsection{Neural Network surrogate}

The autothermal reformer was also replaced by a surrogate block containing a neural network. Hyperparameter tuning was performed to obtain the neural network with the lowest validation loss, quantified with the mean squared error. Hyperparameter ranges are shown in Table 2. The sigmoid and tanh activations were chosen because they are smooth functions, matching our setting of nonlinear continuous optimization. The Adam optimizer was used for training 

\begin{table}[ht]
\centering
\begin{tabular}{ |c|c| } 
 \hline
 \textbf{Hyperparameter} & \textbf{Range/value} \\
 \hline
 Activation function & Sigmoid \& tanh \\ 
 Number of layers & 2 - 5 \\
 Number of neurons & 20 - 35 \\
 Epochs & 500 \\
 \hline
\end{tabular}
\caption{Hyperparameter ranges used to train the neural network.}
\label{table:2}
\end{table}

The training time to run every hyperparameter combination in Table 2 was 550 s. Ultimately, the neural network that best approximates the Gibbs reactor uses the tanh activation function and has 4 hidden layers with 32 neurons each. The time to train this neural network in isolation was 9.5 s. The optimization formulation uses the full space option provided by OMLT, as it was found to converge more reliably than the reduced space formulation. In the full-space formulation, variables and constraints corresponding to interior nodes in the neural network are explicitly included in the optimization problem. While only a single neural network surrogate is considered in this work, a comparison of optimization problems with many different architectures (and equation-based formulations) of embedded neural networks would be an interesting study.

\subsection{Implicit Function}

The theoretical formulation given in Eq. (\ref{eq3}) can be implemented as shown in Eq. (\ref{eq10}). 

\begin{subequations}\label{eq10}
  \begin{align}
    \max &~~ f(x, u) \\
    \text{s.t.} &~~ G(x, u) \leq 0 \\
    &~~ H_{\rm internal}(x, u) = 0 \\
    &~~ H_{\rm linking}(x, y, u) = 0 \\
    &~~ R(x,y,u) = 0
  \end{align}
\end{subequations}

Here, Eq. (\ref{eq10}e) solves for $y$ externally as an implicit function
$y = R_{y}(x,u)$ using the PyNumero interface and the resulting values are communicated back to Eq. (\ref{eq10}d), which is exposed to the NLP solver. In this case, Eq. (\ref{eq10}d) is referred to as a set of linking equality constraints that link the externally obtained outlet variables $y$ to the inlet variables of the reformer recuperator, which is downstream of the Gibbs reactor (See Figure \ref{fig:figure 1}). The dimension of $y$ equals the dimension of
$R$ and $\nabla_y R$ is nonsingular.

The Gibbs reactor is replaced by an implicit function $y = R_{y}(x,u)$ that solves the reactor equations as a parameterized system of equations. The PyNumero interface solves this system at every iteration of the nonlinear optimization solver with a decomposition that partitions variables and equations into block-lower triangular form before solving the resulting blocks independently. The block triangular partition is computed by the approach of Duff and Reid \cite{duff1978implementation} using the Incidence Analysis Pyomo extension \cite{parker2023dulmage}. The PyNumero interface also computes the derivatives required by IPOPT

\begin{figure}[h]
\centering
\resizebox{0.45\columnwidth}{!}{\includegraphics{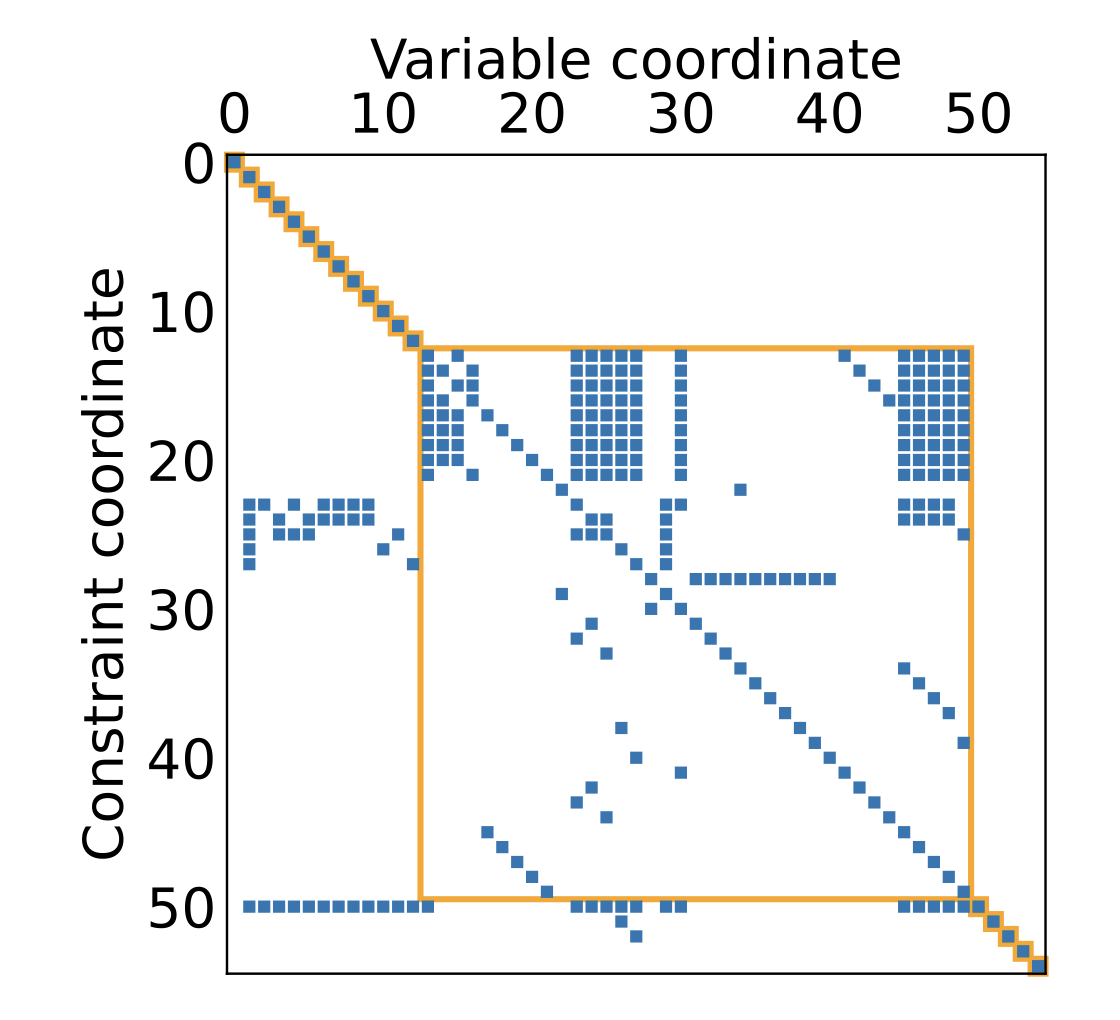}}
\caption{\label{fig:figure 2} Incidence matrix of the square system corresponding to the Gibbs reactor.}
\end{figure}

The system of equations for $R_{y}(x,u)$ decomposes into 55 diagonal blocks, where
18 have dimension 1$\times$1 and one has dimension 37$\times$37.
The block triangular form of this system's incidence matrix is shown in Figure
\ref{fig:figure 2}.

Each independent system of equations is solved with SciPy's {\tt fsolve},
a wrapper around MINPACK's implementation of Powell's hybrid trust region method
\cite{powell1970hybrid}.

\section{Results} 

To compare convergence reliability, solve time and solution quality between the full space, implicit function, and surrogate-based formulations, we perform a parameter sweep varying the inlet natural gas absolute pressure and its conversion in the Gibbs reactor. We attempt to solve the optimization problem using IPOPT \cite{Wchter2005} for every combination of these two parameters for a total of 64 problem instances. These have identical initialization methods and solver options, where the maximum number of iterations are 300 and each must converge to a tolerance of $10^{-7}$. The code used to implement these formulations and reproduce the results can be found at https://github.com/Robbybp/surrogate-vs-implicit. 

\begin{figure}[ht]
\centering
\resizebox{0.95\columnwidth}{!}{\includegraphics{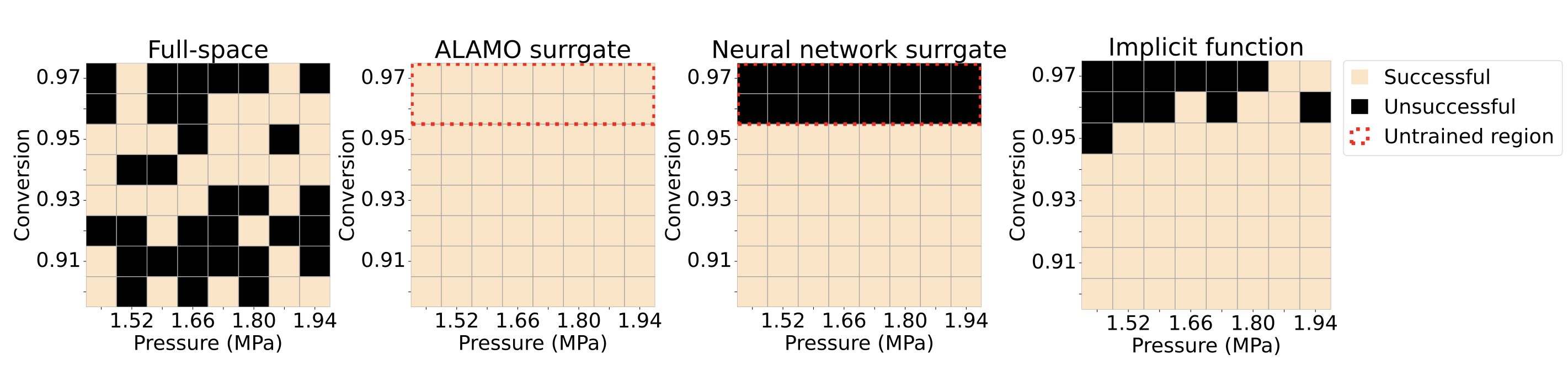}}
\caption{\label{fig:figure 3} Convergence status for each formulation. ``Untrained region'' indicates conversions above 0.95, the upper bound used for surrogate training data. The results indicate that while all three alternative formulations are more reliable than the full-space formulation, the ALAMO surrogate is the most reliable.}
\end{figure}

The convergence status for each formulation is shown in Figure \ref{fig:figure 3}. An unsuccessful run is due to the optimization solver reaching the iteration limit, converging to an infeasible point, or a failure due to repeated function evaluation errors. Here, a function evaluation error may be an error in a scalar-valued function, such as attempting to evaluate the logarithm of a negative number, or a more complicated error such as a failure to solve the square system that defines the implicit function $R_y(x,u)$. A successful run indicates not only that the optimization problem converged, but also that the calculated input variables yield no constraint violations when used to simulate the full-space model.

In this experiment, the full space formulation was able to successfully converge 33 out of 64 instances, the implicit function formulation solved 52, the surrogate-based formulation using ALAMO converged 64 instances, and the surrogate-based formulation using a neural network converged 48. Regarding the full space formulation, many of these failed instances are due to large residuals in calculating energy balances in the autothermal reformer and the reformer recuperator.

The implicit function formulation obtains the same objective and values for the manipulated inputs as the full space formulation but converged 58\% more instances.  The unsuccessful instances at higher conversion are due to function evaluation errors. 

In the failing instances, the implicit function formulation experiences evaluation errors in the Gibbs minimization equations. Given the high conversion, the optimization algorithm calculates a value near zero for the outlet molar composition of propane in the reactor, which is used to calculate the entropy of the ideal gas. Here, one of the terms attempts to calculate the natural logarithm of this composition, causing numerical stability issues, as shown in Eq. (\ref{eq12}). In consequence, the Newton solver fails to calculate $g_{partial,\ch{C3H8}}$, see Eq. (\ref{eq8}), and the outer optimization does not converge. The entropy of the ideal gas is used to calculate the ideal partial Gibbs energy of propane, which is then corrected with a departure function.

\begin{equation}
  s_i^0 = \int_{298.15}^T \frac{A + BT + CT^2 + DT^3}{T} dT + \Delta s_{form}^{298.15} - R\ln\biggl(\frac{P}{P_{ref}}\biggl) - R\ln x_i , ~~i=\ch{C3H8}\label{eq12}
\end{equation}

The ALAMO and neural network formulations introduce only a minor increase in solution error. To evaluate solution quality, we solved a square system for the original flowsheet, where the fixed inlet natural gas, steam molar flow rate, and bypass fraction correspond to the values calculated by the ALAMO and the neural network formulations, separately. Then, we compared the objective value given by these simulations with the objective value calculated by the full space optimization problem. This comparison was done for the 25 instances for which full space, ALAMO, neural network, and implicit formulations all converged successfully. The average relative objective function difference for the ALAMO formulation was 1.9\%, with a maximum difference of 2.1\%. For the neural network formulation, the average relative objective function difference was 0.96\%, with a maximum difference of 2.4\%. 

It is important to note that even though the ALAMO models were trained in the range of conversion of 80 to 95\%, they were able to find accurate solutions with conversions of up to at least 97\%. We note that there is no increase in error when conversions above 95\% are tested.  However, this situation may not be replicated when modeling another process. The neural network surrogate formulation provides virtually the same solution quality as the embedded ALAMO polynomials. Nonetheless, in contrast to those surrogates, it fails to converge when the conversion is 96\% or higher, that is, outside training bounds. One interpretation of these convergence failures is that the neural network surrogate has been over-trained to accurately represent the reactor model in the training region, while the comparatively simple and sparse surrogate generated by ALAMO has resisted this over-training. Additionally, hyperparameter tuning is required to obtain an accurate neural network, which is a computationally costly task.

A breakdown of problem statistics is shown in Table \ref{table:3}, where ``N. Iterations'' refers to the average number of iterations it takes for the optimization solver to converge and ``Std. dev (s)'' refers to the standard deviation of the average solve time from the successful instances. As before, these statistics were computed for the intersection of successful instances among the four formulations.


\begin{table}[ht]
\centering
\resizebox{\textwidth}{!}{
\begin{tabular}{ |c|c|c|c|c|c|c| } 
 \hline
 \textbf{Formulation} & \textbf{N. Iterations} & \textbf{Avg. solve time (s)} & \textbf{Max. solve time (s)} & \textbf{Std. dev. (s)} \\
 \hline
Full Space  & 95 & 1.7 & 5.1 & 0.9 \\ 
Implicit  & 46 & 1.5 & 3.2 & 0.4 \\ 
ALAMO  & 30 & 0.5 & 0.6 & 0.03 \\ 
Neural Network  & 51 & 0.8 & 0.9 & 0.05 \\ 
 \hline
\end{tabular}
}
\caption{Problem statistics for each formulation.}
\label{table:3}
\end{table}

Despite relatively expensive function evaluations that solve a square system at every iteration, the implicit function formulation converges slightly faster than the full-space formulation. This is because the implicit function formulation requires fewer IPOPT iterations to converge the optimization problem. 

The standard deviation of solve time for the full-space formulation shows that different flowsheet parameters have a significant impact on the optimization solver’s ability to converge. Conversely, the implicit and especially surrogate formulations exhibit more uniform solve times and converge for a higher percentage of instances, supporting the idea that they are robust and less sensitive to process parameter values.

Finally, a summary of qualitative results for the four formulations is shown in Table 4. Our qualitative assessment of training time includes time required for data generation and hyperparameter tuning. We note that surrogate training time could be further reduced by using more efficient sampling techniques, such as Latin Hypercube Sampling \cite{kamath}, which would potentially generate a smaller and more representative dataset of the entire experimental region, which might also lead to higher solution accuracies. 

\begin{table}[ht]
\centering
\begin{tabular}{ |c|c|c|c|c| } 
 \hline
  \textbf{Formulation} & \textbf{Solution accuracy} & \textbf{Solve time} & \textbf{Training time} & \textbf{Reliability} \\
 \hline
Full Space & High & Moderate & N/A & Low\\ 
Implicit & High & Moderate & N/A & Moderate\\
ALAMO & Moderate & Low & Moderate & High\\
Neural Network & Moderate & Low & High & Moderate\\ 
 \hline
\end{tabular}
\caption{Qualitative results for each formulation applied to the autothermal reformer
optimization problem.}
\label{table:4}
\end{table}

\section{Conclusions}

We have presented four different formulations for optimization of a chemical process flowsheet using the IDAES modeling framework. The implicit function approach demonstrates an improved convergence reliability in contrast to the full space approach. The ALAMO surrogate formulation is the fastest and most reliable optimization alternative in this study that results in low solution errors for the objective value and the manipulated inputs. Nonetheless, the implicit function formulation may be preferred in cases where the design specifications of a process are not able to tolerate the introduction of errors into the calculations, or where producing enough simulation data to train an accurate surrogate model is computationally prohibitive. For instance, large chemical process flowsheets involving wastewater treatment units or gas scrubber systems, where environmental regulations specify a strict threshold for pollutant compositions in water or gas being released into the atmosphere, may not be appropriate for a surrogate formulation.

In this case study, the ALAMO surrogates were able to effectively approximate the behavior of a Gibbs Reactor. Nevertheless, more complex unit models with polar components, such as multi-component distillation columns, or stripping and absorbing columns, might require neural networks to provide a good approximation of the associated differential-algebraic systems. Finally, the implicit function and surrogate formulations presented in this paper constitute important approaches to optimize a chemical process flowsheet when the classical method fails, and their advantages would be more prominent in the design phase of complex, large scale chemical processes where superstructure optimization (MINLP) might be involved.


\section*{Acknowledgements}
We gratefully acknowledge the support of the U.S. Department of Energy
through the Los Alamos National Laboratory (LANL) LDRD program
and the Center for Nonlinear Studies (CNLS) for this work.

Los Alamos National Laboratory is operated by Triad National Security, LLC,
for the National Nuclear Security Administration of U.S. Department of Energy
(Contract No. 89233218CNA000001).
This work is approved for unlimited release under LA-UR-23-31036.

\pagebreak

\bibliographystyle{ieeetr}
\bibliography{References.bib}

\end{document}